\documentclass{article}
\usepackage{amsmath,amssymb,graphicx}

\newcommand{\nocomma}{}
\newcommand{\tmop}[1]{\ensuremath{\operatorname{#1}}}
\newcommand{\tmtextbf}[1]{{\bfseries{#1}}}
\newcommand{\tmtextit}[1]{{\itshape{#1}}}
\newcommand{\tmtextup}[1]{{\upshape{#1}}}
\newenvironment{tmindent}{\begin{tmparmod}{1.5em}{0pt}{0pt} }{\end{tmparmod}}
\newenvironment{tmparmod}[3]{\begin{list}{}{\setlength{\topsep}{0pt}\setlength{\leftmargin}{#1}\setlength{\rightmargin}{#2}\setlength{\parindent}{#3}\setlength{\listparindent}{\parindent}\setlength{\itemindent}{\parindent}\setlength{\parsep}{\parskip}} \item[]}{\end{list}}
\newenvironment{tmparsep}[1]{\begingroup\setlength{\parskip}{#1}}{\endgroup}
\newtheorem{theorem}{Theorem}

\begin{document}

\title{A Successive Approximation Algorithm for Computing the Divisor
Summatory Function (draft)}\author{Richard Sladkey}\maketitle

\begin{abstract}
  An algorithm is presented to compute isolated values of the divisor
  summatory function in $O \left( n^{1 / 3} \right)$time and $O \left( \log n
  \right)$ space. The algorithm is elementary and uses a geometric approach of
  successive approximation combined with coordinate transformation.
\end{abstract}

\section{Introduction}

Consider the hyperbola from Dirichlet's divisor problem in an $x y$ coordinate
system:
\begin{eqnarray*}
  H \left( x, y \right) & = & x \hspace{0.25em} y\\
  & = & n
\end{eqnarray*}
The number of lattice points under the hyperbola can be thought of as the
number of combinations of positive integers $x$ and $y$ such that their
product is less than or equal to $n$:
\begin{equation}
  T \left( n \right) = \sum_{x, y : xy \leq n} 1
\end{equation}

As such, the hyperbola also represents the divisor summatory function, or the
sum of the number of divisors of all numbers less than or equal to $n$:
\begin{eqnarray}
  \tau \left( x \right) & = & \sum_{d \left| x \right.} 1 = \sum_{x, y : xy =
  n} 1 \\
  T \left( n \right) & = & \sum_{x = 1}^n \tau \left( x \right) \nonumber
\end{eqnarray}
One geometric algorithm is to sum columns of lattice points by choosing an
axis and solving for the variable of the other axis:
\begin{equation}
  T \left( n \right) = \sum^n_{x = 1} \left\lfloor \frac{n}{x} \right\rfloor
\end{equation}
which gives an $O \left( n \right)$ algorithm. By using the symmetry of the
hyperbola (and taking care to avoid double counting) we can do this even more
efficiently:
\begin{equation}
  T \left( n \right) = 2 \sum^{\left\lfloor \sqrt{n} \right\rfloor}_{x = 1}
  \left\lfloor \frac{n}{x} \right\rfloor - \left\lfloor \sqrt{n}
  \right\rfloor^2
\end{equation}
which gives an $O \left( n^{1 / 2} \right)$ algorithm and is in fact the
standard method by which the divisor summatory function is computed. Our goal
is to break this square-root barrier.

In 1903, Vorono\"{\i} in [\ref{bib:Vor03}] made the first significant advance
since Dirichlet on the bound on error term for the divisor problem by
decomposing the hyperbola into a series of non-overlapping triangles
corresponding to tangent lines whose slopes are extended Farey neighbors. We
will use a similar approach but where Vorono\"{\i} produced an exact
expression for the error term and estimated its magnitude, we will instead
produce an algorithm to determine a precise lattice count for an isolated
value of $n$.

\section{Preliminaries}

It will be convenient to parameterize the sum in $T \left( n \right)$ as:
\begin{equation}
  S \left( n, x_{1,} x_2 \right) = \sum^{x_2}_{x = x_1} \left\lfloor
  \frac{n}{x} \right\rfloor
\end{equation}
so that:
\begin{equation}
  T \left( n \right) = S \left( n, 1, n \right) = 2 S \left( n, 1,
  \left\lfloor \sqrt{n} \right\rfloor \right) - \left\lfloor \sqrt{n}
  \right\rfloor^2
\end{equation}
We will also need to count lattice points in triangles. Consider an isosceles
right triangle $\left( 0, 0 \right) \nocomma, \left( i, i \right), \left( i, 0
\right)$, $i$ an integer, excluding points on the bottom gives $1 + 2 + \ldots
+ i$ or:
\[ \Delta \left( i \right) = \frac{i \left( i + 1 \right)}{2} \]
\includegraphics{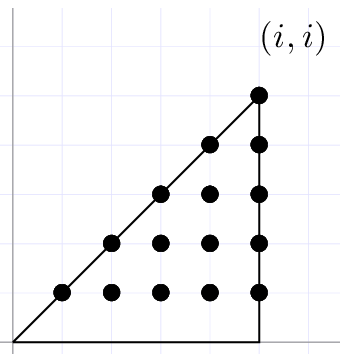}

This formula is also applicable to triangles of the form $\left( 0, 0 \right),
\left( i, ai \right), \left( i, \left( a - 1 \right) i \right)$, $a$ a
positive integer. If we desire to to omit the lattice points on two sides, we
can use $\Delta \left( i - 1 \right)$ instead of $\Delta \left( i \right)$.

\section{Region Processing}

Instead of addressing all of the lattice points, let us for the moment
consider the sub-task of counting the lattice points in a curvilinear
triangular region bounded by two tangent lines and a segment of the hyperbola.
If we can approximate the hyperbola by a series of tangent lines, then the
area below the lines is a simple polygon and can be calculated directly by
decomposing the area into triangles. On the other hand, the region above the
two lines can be handled by chopping off another triangle with a third tangent
line which creates two smaller curvilinear triangular regions.

We will now go about counting the lattice points in such region. We will do
this by first transforming the region into a new coordinate system. This is
very simple conceptually but there are a number of details to take care of in
order to count lattice points accurately and efficiently. First, the tangent
lines are not true tangent lines but are actually shifted to pass through the
nearest lattice points. Because of this, tangent lines need to be ``broken''
on either side of the true tangent point in order to keep them under but close
to the hyperbola. Second, the coordinate transformation turns our simple $xy =
n$ hyperbola into a general quadratic in two variables. Nevertheless, the
recipe at a high level is simply ``tangent, tangent, chop, recurse.''

This figure depicts a typical region in the $x y$ coordinate system:

\includegraphics{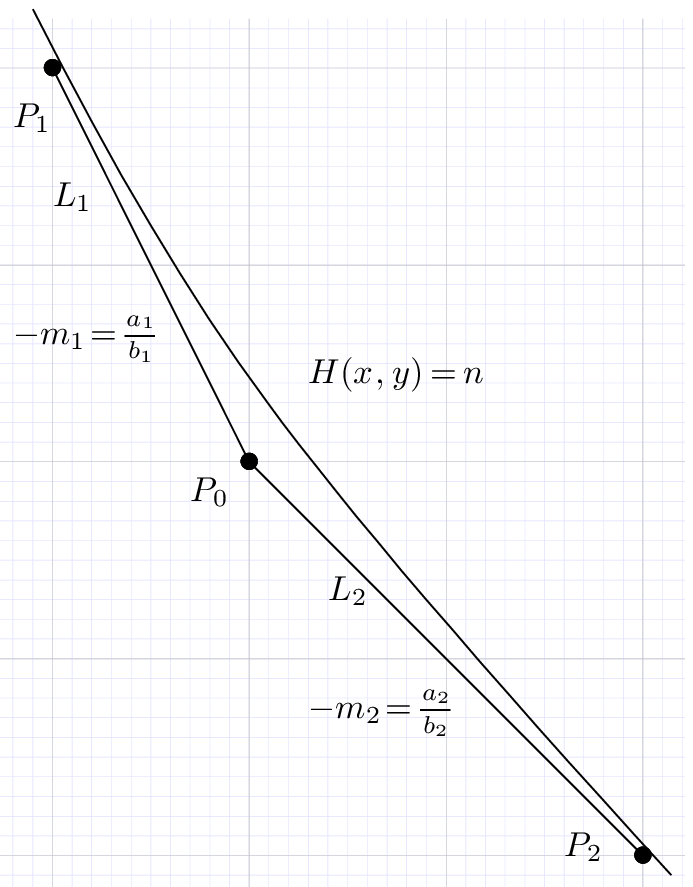}

Define two lines $L_1$ and $L_2$ whose slopes when negated have positive
integral numerators $a_i$ and denominators $b_i$:
\begin{eqnarray}
  - m_1 & = & \frac{a_1}{b_1} \\
  - m_2 & = & \frac{a_2}{b_2} 
\end{eqnarray}
The slopes are chosen to be Farey neighbors so that the determinant is unity:
\begin{equation}
  \left|\begin{array}{cc}
    a_1 & b_1\\
    a_2 & b_2
  \end{array}\right| = a_1  \hspace{0.25em} b_2 - b_1  \hspace{0.25em} a_2 = 1
  \label{eq:det}
\end{equation}
and the slopes are rational numbers which we require to be in lowest terms and
so we can assume $\gcd \left( a_1, b_1 \right) = \gcd \left( a_2, b_2 \right)
= 1$.

Assume further that the lines intersect at the lattice point $P_0$:
\begin{equation}
  \left( x_0, y_0 \right) \label{eq:xy0}
\end{equation}
with $x_0$ and $y_0$ positive integers.

Then the equations for the lines $L_1$ and $L_2$ in point-slope form are:
\begin{eqnarray}
  \frac{y - y_0}{x - x_0} & = & - \frac{a_1}{b_1}  \label{eq:ps1}\\
  \frac{y - y_0}{x - x_0} & = & - \frac{a_2}{b_2}  \label{eq:ps2}
\end{eqnarray}
and converting to standard form:
\begin{eqnarray}
  a_1  \hspace{0.25em} x + b_1  \hspace{0.25em} y & = & x_0  \hspace{0.25em}
  a_1 + y_0  \hspace{0.25em} b_1 \\
  a_2  \hspace{0.25em} x + b_2  \hspace{0.25em} y & = & x_0  \hspace{0.25em}
  a_2 + y_0  \hspace{0.25em} b_2 
\end{eqnarray}
and defining:
\begin{equation}
  c_i = x_0  \hspace{0.25em} a_i + y_0  \hspace{0.25em} b_i
\end{equation}
we have:
\begin{eqnarray}
  a_1  \hspace{0.25em} x + b_1  \hspace{0.25em} y & = & c_1 \\
  a_2  \hspace{0.25em} x + b_2  \hspace{0.25em} y & = & c_2 
\end{eqnarray}
Solving the definitions of $c_1$ and $c_2$ for $x_0$ and $y_0$ give:
\begin{eqnarray}
  x_0 & = & c_1  \hspace{0.25em} b_2 - b_1  \hspace{0.25em} c_2 \\
  y_0 & = & a_1  \hspace{0.25em} c_2 - c_1  \hspace{0.25em} a_2 
\end{eqnarray}
Now observe that the $x y$ lattice points form an alternate lattice relative
to lines $L_1$ and $L_2$:

\includegraphics{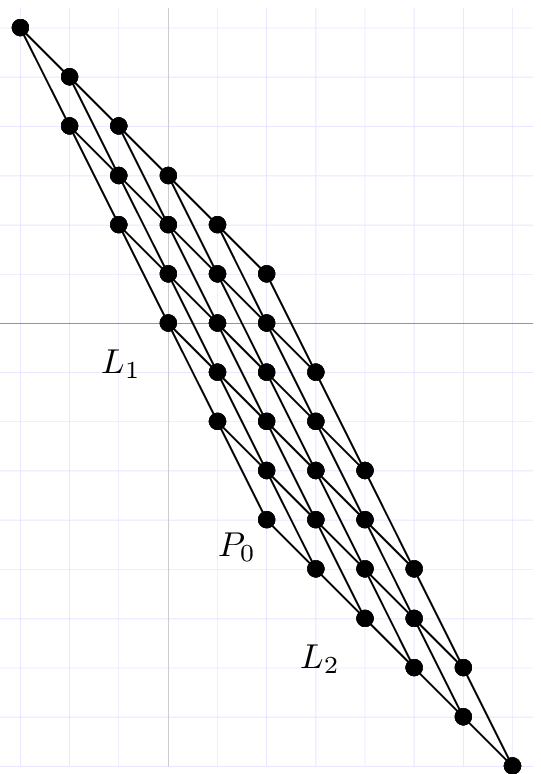}

Define a $u v$ coordinate system with an origin of $P_0$, $L_1$ as the $v$
axis and $L_2$ as the $u$ axis and $u$ and $v$ increasing by one for each
lattice point in the direction of the hyperbola. Then the conversion from the
$u v$ coordinates to $x y$ coordinates is given by:
\begin{eqnarray}
  x & = & x_0 + b_2  \hspace{0.25em} u - b_1  \hspace{0.25em} v 
  \label{eq:uv2xy1}\\
  y & = & y_0 - a_2  \hspace{0.25em} u + a_1  \hspace{0.25em} v 
  \label{eq:uv2xy2}
\end{eqnarray}
Substituting for $x_0$ and $y_0$ and rearranging gives:
\begin{eqnarray}
  x & = & b_2  \hspace{0.25em} \left( u + c_1 \right) - b_1  \hspace{0.25em}
  \left( v + c_2 \right)  \label{eq:uv2xy3}\\
  y & = & a_1  \hspace{0.25em} \left( v + c_2 \right) - a_2  \hspace{0.25em}
  \left( u + c_1 \right)  \label{eq:uv2xy4}
\end{eqnarray}
Solving these equations for $u$ and $v$ and substituting unity for the
determinant provides the inverse conversion from $x y$ coordinates to $u v$
coordinates:
\begin{eqnarray}
  u & = & a_1  \hspace{0.25em} x + b_1  \hspace{0.25em} y - c_1 
  \label{eq:xy2uv1}\\
  v & = & a_2  \hspace{0.25em} x + b_2  \hspace{0.25em} y - c_2 
  \label{eq:xy2uv2}
\end{eqnarray}
Because all quantities are integers, equations (\ref{eq:uv2xy3}),
(\ref{eq:uv2xy4}), (\ref{eq:xy2uv1}), (\ref{eq:xy2uv2}) mean that each $x y$
lattice point corresponds to a $u v$ lattice point and vice versa. As a
result, we can choose to count lattice points in either $x y$ coordinates or
$u v$ coordinates.

Now we are ready to transform the hyperbola into the $u v$ coordinate system
by substituting for $x$ and $y$ in $H \left( x, y \right)$ which gives:
\begin{eqnarray}
  H \left( u, v \right) & = & \left( b_2  \hspace{0.25em} \left( u + c_1
  \right) - b_1  \hspace{0.25em} \left( v + c_2 \right) \right) \left( a_1 
  \hspace{0.25em} \left( v + c_2 \right) - a_2  \hspace{0.25em} \left( u + c_1
  \right) \right) \\
  & = & n \nonumber
\end{eqnarray}
Let us choose a point $P_1$ $\left( 0, h \right)$ on the $v$ axis and a point
$P_2$ $\left( w, 0 \right)$ on the $u$ axis such that:
\begin{eqnarray*}
  H \left( u_h, h \right) & = & n\\
  H \left( w, v_w \right) & = & n\\
  0 \leqslant u_h & < & 1\\
  0 \leqslant v_w & < & 1\\
  - d v / d u \left( u_h \right) & \geq & 0\\
  - d u / d v \left( v_w \right) & \geq & 0
\end{eqnarray*}
or equivalently that the hyperbola is less than one unit away from the nearest
axis at $P_1$ and $P_2$ and that the distance to the hyperbola increases as
you approach the origin.

With these constraints, the hyperbolic segment has the same basic shape as the
full hyperbola: roughly tangent to the axes at the endpoints and strictly
decreasing relative to either axis.

This figure depicts a region in the $u v$ coordinate system:

\includegraphics{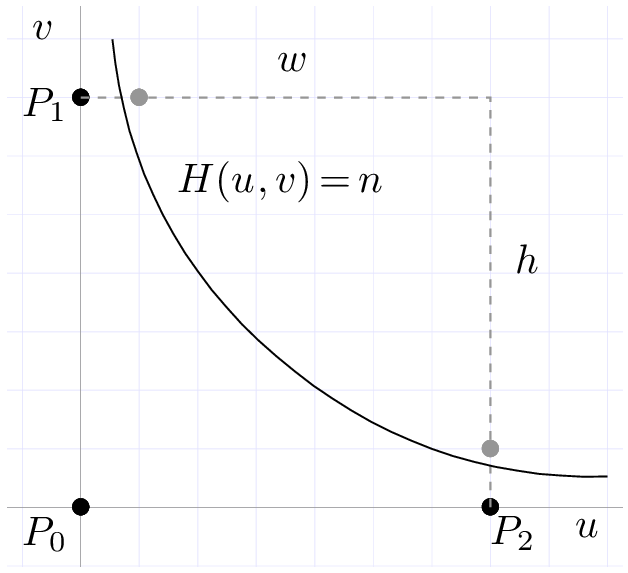}

We can now reformulate the number of lattice points in this region $R_{} $ as
a function of the eight values that define it:
\begin{equation}
  S_{_R} = S_R \left( w, h, a_1, b_1, c_1, a_2, b_2, c_2 \right)
\end{equation}
If $H \left( w, 1 \right) \leq n$, then $v_w \geq 1$ and we can remove the
first lattice row:
\begin{equation}
  S_R = S_R \left( w, h - 1, a_1, b_1, c_1, a_2, b_2, c_2 + 1 \right) + w
\end{equation}
and if $H \left( 1, h \right) \leq n$, then $u_h \geq 1$ and we can remove the
first lattice column:
\begin{equation}
  S_R = S_R \left( w - 1, h, a_1, b_1, c_1 + 1, a_2, b_2, c_2 \right) + h
\end{equation}
so that the conditions are satisified.

At this point we could count lattice points in the region bounded by the $u$
and $v$ axes and $u = w$ and $v = h$ using brute force:
\begin{equation}
  S_R = \sum_{u, v : H \left( u, v \right) \leqslant n} 1
\end{equation}
More efficiently, if we had a formulas for $u$ and $v$ in terms each other, we
could sum columns of lattice points:
\begin{eqnarray}
  & S_W \left( w \right) = \sum^w_{u = 1} \left\lfloor V \left( u \right)
  \right\rfloor & \\
  & S_H \left( h \right) = \sum_{v = 1}^h \left\lfloor U \left( v \right)
  \right\rfloor & 
\end{eqnarray}
using whichever axis has fewer points, keeping in mind that it could be
assymmetric. (Note that these summations are certain not to overcount because
by our conditions $V \left( u \right) < h$ for $0 < u \leq w$ and $U \left( v
\right) < w$ for $0 < v \leq h$.)

And so:
\begin{eqnarray}
  S_R \left( w, h, a_1, b_1, c_1, a_2, b_2, c_2 \right) & = & S_W \\
  & = & S_H \nonumber
\end{eqnarray}
In fact we can derive formulas for $u$ and $v$ in terms of each other by
solving $H \left( u, v \right) = n$ (which when expanded is a general
quadratic in two variables) for $v$ or $u$. The resulting explicit formulas
for $v$ in terms of $u$ and $u$ in terms of $v$ are:
\begin{eqnarray}
  V \left( u \right) & = & \frac{\left( a_1 b_2 + b_1 a_2 \right) 
  \hspace{0.25em} \left( u + c_1 \right) - \sqrt{\left( u + c_1 \right)^2 - 4
  \hspace{0.25em} a_1 b_1 n}}{2 \hspace{0.25em} a_1 b_1} - c_2 \\
  U \left( v \right) & = & \frac{\left( a_1 b_2 + b_1 a_2 \right) 
  \hspace{0.25em} \left( v + c_2 \right) - \sqrt{\left( v + c_2 \right)^2 - 4
  \hspace{0.25em} a_2 b_2 n}}{2 \hspace{0.25em} a_2 b_2} - c_1 
\end{eqnarray}
(Note exchanging $u$ for $v$ results in the same formula with subscripts 1 and
2 exchanged.)

As a result we can compute the number of lattice points within the region
using a method similar to the method usually used for the hyperbola as a
whole. Our goal, however, it to subdivide the region into two smaller regions
and process them recursively, only using manual counting at our discretion. To
do so we need to remove an isosceles right triangle in the lower-left corner
and what will be left are two sub-regions in the upper-left and lower-right.

This figure shows the right triangle and the two sub-regions:

\includegraphics{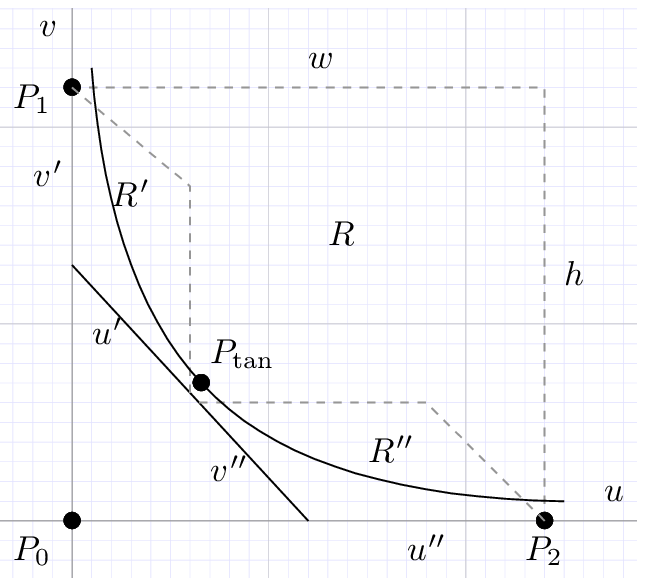}

A diagonal with slope -1 in the $u v$ coordinate system has a slope in the $x
y$ coordinate system that is the mediant of the slopes of lines $L_1$ and
$L_2$:
\begin{eqnarray}
  - m_3 & = & \frac{a_1 + a_2}{b_1 + b_2} 
\end{eqnarray}
So let us define:
\begin{eqnarray}
  a_3 & = & a_1 + a_2 \\
  b_3 & = & b_1 + b_2 
\end{eqnarray}
Then differentiating $H \left( u, v \right) = n$ with respect to $u$ and
setting $d v / d u = - 1$ gives:
\begin{equation}
  \left( a_1  \hspace{0.25em} b_2 + b_1  \hspace{0.25em} a_2 + 2
  \hspace{0.25em} a_2  \hspace{0.25em} b_2 \right)  \hspace{0.25em} \left( u +
  c_1 \right) = \left( a_1  \hspace{0.25em} b_2 + b_1  \hspace{0.25em} a_2 + 2
  \hspace{0.25em} a_1  \hspace{0.25em} b_1 \right)  \hspace{0.25em} \left( v +
  c_2 \right)
\end{equation}
and the intersection of this line with $H \left( u, v \right) = n$ gives the
point $P_{\tan}$ on the hyperbola where the slope is equal to -1:
\begin{eqnarray}
  u_{\tan} & = & \left( a_1  \hspace{0.25em} b_2 + b_1  \hspace{0.25em} a_2 +
  2 \hspace{0.25em} a_1  \hspace{0.25em} b_1 \right)  \hspace{0.25em}
  \sqrt{\frac{n}{a_3  \hspace{0.25em} b_3}} - c_1 \\
  v_{\tan} & = & \left( a_1  \hspace{0.25em} b_2 + b_1  \hspace{0.25em} a_2 +
  2 \hspace{0.25em} a_2  \hspace{0.25em} b_2 \right)  \hspace{0.25em}
  \sqrt{\frac{n}{a_3  \hspace{0.25em} b_3}} - c_2 
\end{eqnarray}
The equation of a line through this intersection and tangent to the hyperbola
is then $u + v = u_{\tan} + v_{\tan}$ which simplifies to:
\begin{equation}
  u + v = 2 \hspace{0.25em} \sqrt{a_3  \hspace{0.25em} b_3  \hspace{0.25em} n}
  - c_1 - c_2
\end{equation}
Next we need to find the pair of lattice points $P_4$ $\left( u_{4,} v_4
\right)$ and $P_5$ $\left( u_5, v_5 \right)$ such that:
\begin{eqnarray*}
  u_4 & > & 0\\
  u_5 & = & u_4 + 1\\
  - d v / d u \left( u_4 \right) & \geq & 1\\
  - d v / d u \left( u_5 \right) & < & 1\\
  v_4 & = & \left\lfloor V \left( u_4 \right) \right\rfloor\\
  v_5 & = & \left\lfloor V \left( u_5 \right) \right\rfloor
\end{eqnarray*}
The derivative conditions ensure that the diagonal rays with slope $- 1$
pointing outward from $P_4$ and $P_5$ do not intersect the hyperbola. Setting
$u_4 = \left\lfloor u_{\tan} \right\rfloor$ will satisfy the conditions as
long as $u_4 \neq 0$.

Let the point at which the ray from $P_4$ intersects the $v$ axis be $P_6
\left( 0, v_6 \right)$ and the point at which the ray from $P_5$ intersects
the $u$ axis be $P_7 \left( u_7, 0 \right)$. Then:
\begin{eqnarray}
  v_6 & = & u_4 + v_4 \\
  u_7 & = & u_5 + v_5 
\end{eqnarray}
A diagram of all the points defined so far:

\includegraphics{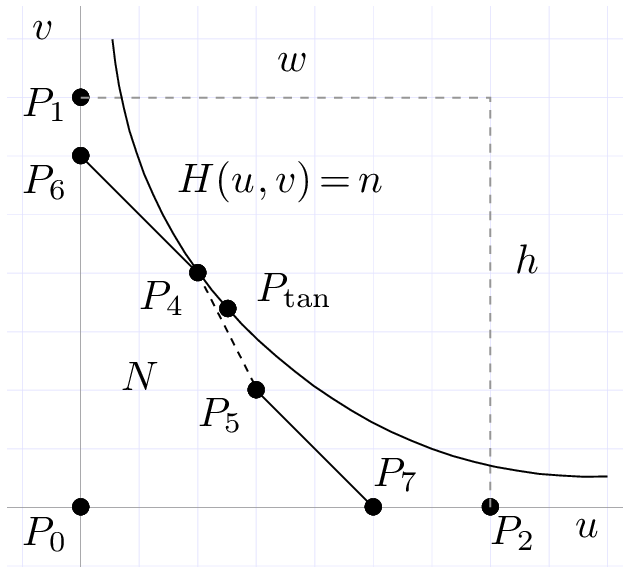}

Then the number of lattice points above the axes and inside the polygon $N$
defined by points $P_0, P_6, P_4, P_5, P_7$ is
\[ S_N = \Delta \left( v_6 - 1 \right) - \Delta \left( v_6 - u_5 \right) +
   \Delta \left( u_7 - u_5 \right) \]
or
\begin{equation}
  S_N = \left\{ \begin{array}{ll}
    \Delta \left( v_6 - 1 \right) + u_4 & \tmop{if} v_6 < u_7\\
    \Delta \left( v_6 - 1 \right) & \tmop{if} v_6 = u_7\\
    \Delta \left( u_7 - 1 \right) + v_5 & \tmop{if} v_6 > u_7
  \end{array} \right.
\end{equation}
because counting on reverse lattice diagonals starting at the origin we sum $1
+ 2 + \ldots + \left( \min \left( v_6, u_7 \right) - 1 \right)$ plus a partial
diagonal if the polygon is not a triangle.

Using the properties of Farey fractions observe that:
\begin{eqnarray*}
  \left|\begin{array}{cc}
    a_1 & b_1\\
    a_3 & b_3
  \end{array}\right| & = & a_1  \left( b_1 + b_2 \right)_{} - b_1  \left( a_1
  + a_2 \right) = a_1 b_2 - b_1 a_2 = \left|\begin{array}{cc}
    a_1 & b_1\\
    a_2 & b_2
  \end{array}\right| = 1\\
  \left|\begin{array}{cc}
    a_3 & b_3\\
    a_2 & b_2
  \end{array}\right| & = &  \left( a_1 + a_2 \right)_{} b_2 - \left( b_1 + b_2
  \right) a_2 = a_1 b_2 - b_1 a_2 = \left|\begin{array}{cc}
    a_1 & b_1\\
    a_2 & b_2
  \end{array}\right| = 1
\end{eqnarray*}
so that $m_1$ and $m_3$ are also Farey neighbors and likewise for $m_3$ and
$m_2$.

So we can define region $R'$ to be the sub-region with $P'_1 = P_1, P'_0 =
P_6, P'_2 = P_4$ and the region $R''$ to be the sub-region with $P''_1 = P_5,
P''_0 = P_7, P_2'' = P_2$ and then the number of lattice points in the entire
region is $S_R = S_N + S_{R'} + S_{R''}$ or
\begin{eqnarray}
  S_R \left( w, h, a_1, b_1, c_1, a_2, b_2, c_2 \right) & = & S_N \\
  & + & S_R \left( u_4, h - v_6, a_1, b_1, c_1, a_3, b_3, c_1 + c_2 + v_6
  \right) \nonumber\\
  & + & S_R \left( w - u_7, v_5, a_3, b_3, c_1 + c_2 + u_7, a_2, b_2, c_2
  \right) . \nonumber
\end{eqnarray}
This recursive formula for the sum of the lattice points in a region in terms
of the lattice points in its sub-regions allows us to use a divide and conquer
approach to counting lattice points under the hyperbola.

\section{Top Level Processing}

Now let us return to the hyperbola as a whole. It should be clear that it is
easy in $x y$ coordinates to calculate $y$ in terms of $x$ by solving $H
\left( x, y \right) = n$ for $y$:
\begin{equation}
  Y \left( x \right) = \frac{n}{x}
\end{equation}
We know that we only need to sum lattice points under the hyperbola up to
$\left\lfloor \sqrt{n} \right\rfloor$. The point $\sqrt{n}$ is in fact at the
$x = y$ axis of symmetry and so the slope at that point is exactly $- 1$. The
next integral slope occurs at $- 2$, so our first (and largest) region occurs
between slopes $- m_1 = 2$ and $- m_2 = 1$. By processing adjacent integral
slopes we will start in the middle and work our way back towards the origin.

However, we cannot use the region method for the whole hyperbola because
regions become smaller and smaller and eventually a region has a size $w + h
\leq 1$. We can find the point where this occurs by taking the second
derivative of $Y \left( x \right)$ with respect to $x$ and setting it to
unity. In other words, the point on the hyperbola where the rate of change in
the slope exceeds one per lattice column, which is:
\begin{equation}
  x = \sqrt[3]{2 n} = 2^{1 / 3} n^{1 / 3} \approx 1.26 n^{1 / 3}
\end{equation}
As a result there is no benefit in region processing the first $O \left( n^{1
/ 3} \right)$ lattice columns so we resort to the simple method to sum the
lattice columns less than $x_{\min}$:
\begin{eqnarray}
  x_{\min} & = & C_1  \left\lceil \sqrt[3]{2 n} \right\rceil \\
  x_{\max} & = & \left\lfloor \sqrt{n} \right\rfloor \nonumber\\
  y_{\min} & = & \left\lfloor Y \left( x_{\max} \right) \right\rfloor
  \nonumber\\
  S_1 & = & S \left( n, 1, x_{\min} - 1 \right) \nonumber
\end{eqnarray}
where $C_1 \geq 1$ is a constant to be chosen later.

Next we need to account for the all the points on or below the first line
which is a rectangle and a triangle:
\[ S_2 = \left( x_{\max} - x_{\min} + 1 \right) y_{\min} + \Delta \left(
   x_{\max} - x_{\min} \right) \]
Because all slopes in this section of the algorithm are whole integers, we
have:
\begin{eqnarray*}
  a_i & = & - m_i\\
  b_i & = & 1
\end{eqnarray*}
Assume that we have point $P_2 \nocomma$ and value $a_2$ from the previous
iteration. For the first iteration we will have:
\begin{eqnarray*}
  x_2 & = & x_{\max}\\
  y_2 & = & y_{\min}\\
  a_2 & = & 1
\end{eqnarray*}
For all iterations:
\[ a_1 = a_2 + 1 \]
The $x$ coordinate of the point on the hyperbola where the slope is equal to
$m_1$ can be found by taking the derivative of $Y \left( x \right)$ with
respect to $x$, setting $d y / d x = m_1$, and then solving for $x$:
\begin{equation}
  x_{\tan} = \sqrt{\frac{n}{a_1}}
\end{equation}
Similar to processing a region (but now in $x y$ coordinates), we now need two
lattice points $P_4  \left( x_{4,} y_4 \right)$ and $P_5  \left( x_5, y_5
\right)$ such that:
\begin{eqnarray*}
  x_4 & > & x_{\min}\\
  x_5 & = & x_4 + 1\\
  - d y / d x \left( x_4 \right) & \geq & a_1\\
  - d y / d x \left( x_5 \right) & < & a_1\\
  y_4 & = & \left\lfloor Y \left( x_4 \right) \right\rfloor\\
  y_5 & = & \left\lfloor Y \left( x_5 \right) \right\rfloor
\end{eqnarray*}
To meet these conditions we can set $x_4 = \left\lfloor x_{\tan}
\right\rfloor$ unless $x_4 \leq x_{\min}$ in which case we can manually count
the lattice columns between $x_{\min}$ and $x_2$ and cease iterating. If so,
the remaining columns can be computed as:
\[ S_3 = \sum^{x_{2 - 1}}_{x = x_{\min}} \left\lfloor \frac{n}{x}
   \right\rfloor - \left( a_2  \left( x_2 - x \right) + y_2 \right) \]
which is the number of lattice points below the hyperbola and above line $L_2$
over the interval $\left[ x_{\min}, x_2 \right)$.

Now take line $L_2$ with slope $- a_2$ passing through $P_2$, lines $L_4$ and
$L_5$ with slopes $- a_1$ and passing through $P_4$ and $P_5$ and then find
the point $P_6$ where $L_4$ intersects $x = x_{\min}$ and the point $P_0$
where $L_5$ intersects $L_2$ and the point $P_7$ where $L_2$ intersects $x =
x_{\min}$ and denote by $c_i$ the $y$ intercept of line $L_i$.

\includegraphics{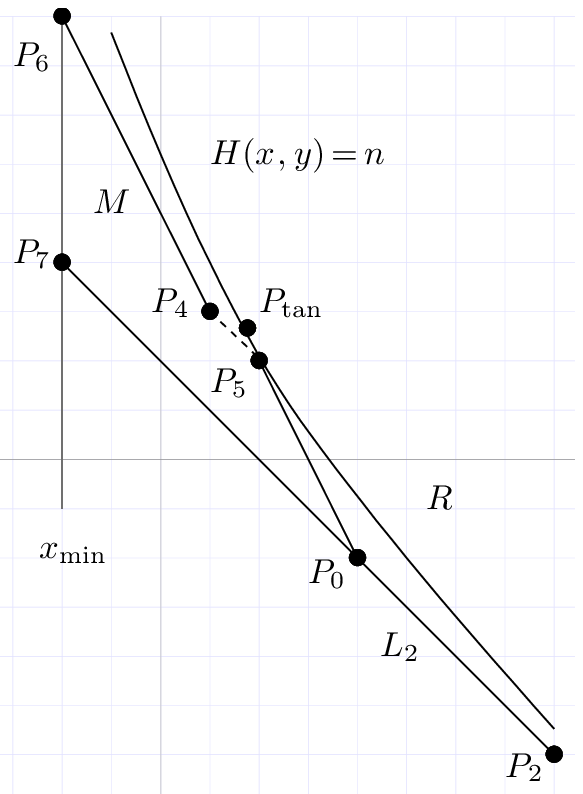}

Now add up the lattice points in the polygon $M$ defined by the points $P_0,
P_7, P_6, P_4, P_5$ but above $L_2$ by adding the whole triangle corresponding
to $L_4$, subtracting the portion of it to the right of $P_4$, and then adding
back the triangle corresponding to $L_5$ stating at $P_5$:
\begin{equation}
  S_M = \Delta \left( c_4 - c_2 - x_{\min} \right) - \Delta \left( c_4 - c_2 -
  x_5 \right) + \Delta \left( c_5 - c_2 - x_5 \right)
\end{equation}
where if $L_4$ is coincident with $L_5$, the second two terms cancel each
other out.

Then choosing $P_1 = P_5$ (together with $P_0$ and $P_2$) and calculating the
necessary quantities we have a region $R$ and can now count lattice points
using region processing:
\begin{equation}
  S_R = S_R \left( a_1 x_2 + y_2 - c_5, a_2 x_5 + y_5 - c_2, a_1, 1, c_5, a_2,
  1, c_2 \right)
\end{equation}
so the total sum for this iteration is:
\[ S_A \left( a_1 \right) = S_M + S_R \]
Then we may advance to the next region by setting:
\begin{eqnarray*}
  x'_2 & = & x_4\\
  y'_2 & = & y_4\\
  a'_2 & = & a_1
\end{eqnarray*}
Summing all interations gives
\[ S_4 = \sum^{a_{\max}}_{a = 2} S_A \left( a \right) . \]
Finally, the total number of lattice points under the hyperbola from $1$ to
$x_{\max}$ is
\begin{equation}
  S_T = S \left( 1, x_{\max} \right) = S_1 + S_2 + S_3 + S_4
\end{equation}
and therefore the final computation of the divisor summatory function is given
by
\begin{equation}
  T \left( n \right) = 2 S_T - \left\lfloor \sqrt{n} \right\rfloor .^2
\end{equation}

\section{Division-Free Counting}

Since we calculate $S_1$ using the traditional method and since the
computation will consist entirely of $S_1$ when $n < 4 C_1^6$, it is
beneficial to have a faster method of performing this step, albeit by a
constant factor. Denote by $l = \left\lceil \log_2 \left( n \right)
\right\rceil$ the number of bits needed to represent $n$. We can avoid an
$l$-bit division in most iterations by using a Bresenham-style calculation
(see [\ref{bib:Bres77}]) and working backwards while computing an estimate of
the result of the division based on the previous iteration.

Define $\beta \left( x \right) = \left\lfloor Y \left( x \right)
\right\rfloor$, the finite difference $\delta_1 \left( x \right) = \beta
\left( x \right) - \beta \left( x + 1 \right)$, and the second-order finite
difference $\delta_2 \left( x \right) = \delta_1 \left( x \right) - \delta_1
\left( x + 1 \right)$. To check whether the value is correct we also need to
keep track of the error. So defining the error $\varepsilon \left( x \right) =
n - x \beta \left( x \right) = n - x \left\lfloor n / x \right\rfloor = n
\tmop{mod} x$ gives
\begin{eqnarray*}
  \varepsilon \left( x \right) - \varepsilon \left( x + 1 \right) & = & \left(
  x + a \right) \beta \left( x + 1 \right) - x \beta \left( x \right)\\
  & = & \left( x + 1 \right) \beta \left( x + 1 \right) - x \left( \beta
  \left( x + 1 \right) + \delta_1 \left( x + 1 \right) + \delta_2 \left( x
  \right) \right)\\
  & = & \beta \left( x + 1 \right) - x \delta_1 \left( x + 1 \right) - x
  \delta_2 \left( x \right)
\end{eqnarray*}
Introducing the intermediate quantity $\gamma \left( x \right) = \beta \left(
x \right) - \left( x - 1 \right) \delta_1 \left( x \right)$ and
$\hat{\varepsilon} \left( x \right)$ as the estimate of the error assuming
$\delta_2 \left( x \right)$=0 then
\begin{eqnarray*}
  \hat{\varepsilon} \left( x \right) & = & \varepsilon \left( x + 1 \right) +
  \gamma \left( x + 1 \right)\\
  \delta_2 \left( x \right) & = & \left\lfloor \frac{\hat{\varepsilon} \left(
  x \right)}{x} \right\rfloor\\
  \delta_1 \left( x \right) & = & \delta_1 \left( x + 1 \right) + \delta_2
  \left( x \right)\\
  \varepsilon \left( x \right) & = & \hat{\varepsilon} \left( x \right) - x
  \delta_2 \left( x \right)\\
  \gamma \left( x \right) & = & \gamma \left( x + 1 \right) + 2 \delta_1
  \left( x \right) - x \delta_2 \left( x \right)\\
  \beta \left( x \right) & = & \beta \left( x + 1 \right) + \delta_1 \left( x
  \right) .
\end{eqnarray*}
Over the range $x_1 \leq x \leq x_2$ these integer quantites are bounded in
size by $x, \varepsilon \left( x \right), \left| \hat{\varepsilon} \left( x
\right) \right| \leq x_2, \left| \gamma \left( x \right) \right| \leq \max
\left( 2 n / \left( x_1^2 + x_1 \right), x_2 \right)$, $\beta \left( x \right)
\leq n / x_1$, $\delta_1 \left( x \right) \leq n / x_1^2$+1,$\left| \delta_2
\left( x \right) \right| \leq 2 n / x_1^3 + 2$.

For $\sqrt[3]{2 n} < x \leq \sqrt{n}$, $\delta_2 \left( x \right) \in \left\{
- 1, 0, 1, 2 \right\}$ and so
\[ \left\lfloor \frac{\hat{\varepsilon} \left( x \right)}{x} \right\rfloor =
   \left\{ \begin{array}{rl}
     2 & \tmop{if} \hat{\varepsilon} \left( x \right) \geq 2 x ;\\
     1 & \tmop{if} x \leq \hat{\varepsilon} \left( x \right) < 2 x ;\\
     - 1 & \tmop{if} \hat{\varepsilon} \left( x \right) < 0 ;\\
     0 & \tmop{otherwise} ;
   \end{array} \right. \]
and thus $\beta \left( x \right), \gamma \left( x \right), \delta_1 \left( x
\right), \varepsilon \left( x \right)$ can be computed from $\beta \left( x +
1 \right), \gamma \left( x + 1 \right), \delta_1 \left( x + 1 \right),
\varepsilon \left( x + 1 \right)$ using only addition and subtraction of
$\frac{1}{2} l$-bit quantities except $\beta \left( x \right)$ which is
$\frac{2}{3} l$ bits. Note that $\hat{\varepsilon} \left( x \right) \geq 2 x$
is very rare over this range and if $\hat{\varepsilon} \left( x \right) \geq 3
x$, it means that $x < \sqrt[3]{2 n} .$ For $n^{1 / 6} \leq x \leq \sqrt[3]{2
n}$ we can add the modest division $\left\lfloor \hat{\varepsilon} \left( x
\right) / x \right\rfloor$ between two $\frac{1}{3} l$-bit values, $\gamma
\left( x \right)$ and $\delta_1 \left( x \right)$ grow to $\frac{2}{3} l$ bits
and $\beta \left( x \right)$ grows to $\frac{5}{6} l$ bits. \ For $x < n^{1 /
6}$ we can sum using ordinary division.

\section{Algorithms}

In this section we present a series of algorithms based on the previous
sections. \ The short-hand notation $F \left( x \right) : \tmop{expression}$
signifies a functional value that remains unevaluated until referenced.

The first algorithm is a straightforward version of the basic successive
approximation method. A literal implementation based on this description will
offer many opportunities for optimization. Various formulas have been slightly
modified so that the entire algorithm can be implemented using only unsigned
multi-precision integer arithmetic. The operations required are addition,
subtraction, multiplication, floor division, floor square root, ceiling square
root, and ceiling cube root. If any of the root operations are not available,
they may be implemented using Newton's method.

{\noindent}\begin{tmparsep}{0em}
  \tmtextbf{Algorithm \tmtextup{1}}{\smallskip}
  
  \begin{tmindent}
    Inputs: $n \geq 0, C_1 \approx 10, C_2 \approx 10$

    $\Delta \left( i \right) : i \left( i + 1 \right) / 2$
    
    $S_1 \left( \right) : \sum_{x = 1}^{x < x_{\min}} \left\lfloor n / x
    \right\rfloor$
    
    $S_2 \left( \right) : \left( x_{\max} - x_{\min} + 1 \right) y_{\min} +
    \Delta \left( x_{\max} - x_{\min} \right)$
    
    $S_3 \left( \right) : \sum_{x = x_{\min}^{}}^{x < x_2} \left\lfloor n / x
    \right\rfloor - \left( a_2  \left( x_2 - x \right) + y_2 \right)$
    
    $S_M \left( \right) : \Delta \left( c_4 - c_2 - x_{\min} \right) - \Delta
    \left( c_4 - c_2 - x_5 \right) + \Delta \left( c_5 - c_2 - x_5 \right)$

    $x_{\max} \leftarrow \left\lfloor \sqrt{n} \right\rfloor, y_{\min}
    \leftarrow \left\lfloor n / x_{\max} \right\rfloor, x_{\min} \leftarrow
    \min \left(  \left\lceil C_1  \sqrt[3]{2 n} \right\rceil, x_{\max}
    \right)$
    
    $s \leftarrow 0, a_2 \leftarrow 1, x_2 \leftarrow x_{\max}, y_2 \leftarrow
    y_{\min}, c_2 \leftarrow a_2 x_2 + y_2$
    
    \tmtextbf{loop}
    \begin{tmindent}
      $a_1 \leftarrow a_2 + 1$
      
      $x_4 \leftarrow \left\lfloor \sqrt{\left\lfloor n / a_1 \right\rfloor}
      \right\rfloor, y_4 \leftarrow \left\lfloor n / x_4 \right\rfloor, c_4
      \leftarrow a_1 x_4 + y_4$
      
      $x_5 \leftarrow x_4 + 1, y_5 \leftarrow \left\lfloor n / x_5
      \right\rfloor, c_5 \leftarrow a_1 x_5 + y_5$
      
      \tmtextbf{if} $x_4 < x_{\min}$ \tmtextbf{then} \tmtextbf{exit}
      \tmtextbf{loop} \tmtextbf{end} \tmtextbf{if}
      
      $s \leftarrow s + S_M \left( \right) + S_R \left( a_1 x_2 + y_2 - c_5,
      a_2 x_5 + y_5 - c_2, a_1, 1, c_5, a_2, 1, c_2 \right)$
      
      $a_2 \leftarrow a_1, x_2 \leftarrow x_4, y_2 \leftarrow y_4, c_2
      \leftarrow c_4$
    \end{tmindent}
    \tmtextbf{end} \tmtextbf{loop}
    
    $s \leftarrow s + S_1 \left( \right) + S_2 \left( \right) + S_3 \left(
    \right)$
    
    \tmtextbf{return} $2 s - x_{\max}^2$

    \tmtextbf{function} $S_R \left( w, h, a_1, b_1, c_1, a_2, b_2, c_2
    \right)$
    \begin{tmindent}
      $\Delta \left( i \right) : i \left( i + 1 \right) / 2$
      
      $H \left( u, v \right) : \left( b_2  \hspace{0.25em} \left( u + c_1
      \right) - b_1  \hspace{0.25em} \left( v + c_2 \right) \right)  \left(
      a_1  \hspace{0.25em} \left( v + c_2 \right) - a_2  \hspace{0.25em}
      \left( u + c_1 \right) \right)$
      
      $U_{\tan} \left( \right) : \left\lfloor \sqrt{\left\lfloor \left( a_1 
      \hspace{0.25em} b_2 + b_1  \hspace{0.25em} a_2 + 2 \hspace{0.25em} a_1 
      \hspace{0.25em} b_1 \right)^2 n / \left( a_3  \hspace{0.25em} b_3
      \right) \right\rfloor} \right\rfloor - c_1$
      
      $V_{\tmop{floor}} \left( u \right) : \left\lfloor \left( \left( a_1 
      \hspace{0.25em} b_2 + b_1  \hspace{0.25em} a_2 \right)  \hspace{0.25em}
      \left( u + c_1 \right) - \left\lceil \sqrt{\left( u + c_1 \right)^2 - 4
      \hspace{0.25em} a_1  \hspace{0.25em} b_1  \hspace{0.25em} n}
      \right\rceil \right) / \left( 2 a_1 b_1 \right) \right\rfloor - c_2$
      
      $U_{\tmop{floor}} \left( v \right) : \left\lfloor \left( \left( a_1 
      \hspace{0.25em} b_2 + b_1  \hspace{0.25em} a_2 \right)  \hspace{0.25em}
      \left( v + c_2 \right) - \left\lceil \sqrt{\left( v + c_2 \right)^2 - 4
      \hspace{0.25em} a_2  \hspace{0.25em} b_2  \hspace{0.25em} n}
      \right\rceil \right) / \left( 2 a_2 b_2 \right) \right\rfloor - c_2$
      
      $S_W \left( \right) : \sum_{u = 1}^{u < w} V_{\tmop{floor}} \left( u
      \right)$
      
      $S_H \left( \right) : \sum_{v = 1}^{v < h} U_{\tmop{floor}} \left( v
      \right)$
      
      $S_N \left( \right) : \Delta \left( v_6 - 1 \right) - \Delta \left( v_6
      - u_5 \right) + \Delta \left( u_7 - u_5 \right)$

      $s \leftarrow 0, a_3 \leftarrow a_1 + a_2, b_3 \leftarrow b_1 + b_2$
      
      \tmtextbf{if} $h > 0 \wedge H \left( w, 1 \right) \leq n$
      \tmtextbf{then} $s \leftarrow s + w, c_2 \leftarrow c_2 + 1, h
      \leftarrow h - 1$ \tmtextbf{end} i\tmtextbf{f}
      
      \tmtextbf{if} $w > 0 \wedge H \left( 1, h \right) \leq n$
      \tmtextbf{then} $s \leftarrow s + h, c_1 \leftarrow c_1 + 1, w
      \leftarrow w - 1$ \tmtextbf{end} i\tmtextbf{f}
      
      \tmtextbf{if} $w \leq C_2$ \tmtextbf{then} \tmtextbf{return} $s + S_W
      \left( \right)$ \tmtextbf{end} \tmtextbf{if}
      
      \tmtextbf{if} $h \leq C_2$ \tmtextbf{then} \tmtextbf{return} $s + S_H
      \left( \right)$ \tmtextbf{end} \tmtextbf{if}
      
      $u_4 \leftarrow U_{\tan} \left( \right), v_4 \leftarrow V_{\tmop{floor}}
      \left( u_4 \right), u_5 \leftarrow u_4 + 1, v_5 \leftarrow
      V_{\tmop{floor}} \left( u_5 \right)$
      
      $v_6 \leftarrow u_4 + v_4, u_7 \leftarrow u_6 + v_6$
      
      $s \leftarrow s + S_N \left( \right)$
      
      $s \leftarrow s + S_R \left( u_4, h - v_6, a_1, b_1, c_1, a_3, b_3, c_1
      + c_2 + v_6 \right)$
      
      $s \leftarrow s + S_R \left( w - u_7, v_5, a_3, b_3, c_1 + c_2 + u_7,
      a_2, b_2, c_2 \right)$
      
      \tmtextbf{return} s
    \end{tmindent}
    \tmtextbf{end} \tmtextbf{function}
  \end{tmindent}
\end{tmparsep}{\hspace*{\fill}}{\medskip}

\label{algorithm1}

The next algorithm gives a flavor for the optimizations that are available. \
It computes the manual summation of a small region over $u$ or $v$ using a
handful of additions, one square root and one division per lattice column. A
similar technique can be used to compute $V_{\tmop{floor}}$ for the adjacent
values $u_4$ and $u_5$. Making this portion of the computation faster favors
larger values of $C_2$, the cutoff for small regions. An analogy is that this
step is faster for small regions in the same way that an insertion sort is
faster than a quicksort for small arrays and the break even point can be
determined experimentally.

{\noindent}\begin{tmparsep}{0em}
  \tmtextbf{Algorithm \tmtextup{2}}{\smallskip}
  
  \begin{tmindent}
    $S_W \left( \right) : S_I \left( w, c_1, c_2, a_1 b_2 + b_1 a_2, 2 a_1 b_1
    \right)$
    
    $S_H \left( \right) : S_I \left( h, c_2, c_1, a_1 b_2 + b_1 a_{2,} 2 a_2
    b_2 \right)$

    \tmtextbf{function} $S_I \left( i_{\max}, p_1, p_2, q, r \right)$
    \begin{tmindent}
      $s \leftarrow 0, A \leftarrow p_1^2 - 2 rn, B \leftarrow p_1 q, C
      \leftarrow 2 p_1 - 1$
      
      \tmtextbf{for} $i = 1, \ldots, i_{\max} - 1$ \tmtextbf{do}
      \begin{tmindent}
        $C \leftarrow C + 2, A \leftarrow A + C, B \leftarrow B + q$
      \end{tmindent}
      \begin{tmindent}
        $s \leftarrow s + \left\lfloor \left( B - \left\lceil \sqrt{A}
        \right\rceil \right) / r \right\rfloor$
      \end{tmindent}
      \tmtextbf{end} \tmtextbf{for}
      
      \tmtextbf{return} $s - \left( i_{\max} - 1 \right) p_2$
    \end{tmindent}
    \tmtextbf{end} \tmtextbf{function}
  \end{tmindent}
\end{tmparsep}{\hspace*{\fill}}{\medskip}

The next algorithm formalizes the steps of the division-free counting method
which can be used for the summation $S_1$. \ Whether this is actually faster
depends on many things but for example if $n < 2^{94}$, then $\beta, \delta,
\left| \gamma \right|, \left| \varepsilon \right| < 2^{63}$ for $2^{32} < x <
2^{47}$ and if signed 64-bit addition is a single-cycle operation, then a
computation of $\beta$ using this method is about ten cycles vs. say a hundred
\ cycles for a single multi-precision division.

{\noindent}\begin{tmparsep}{0em}
  \tmtextbf{Algorithm \tmtextup{3}}{\smallskip}
  
  \begin{tmindent}
    $S_1 \left( \right) : S_Q \left( 1, x_{\min} - 1 \right)$

    \tmtextbf{function} $S_Q \left( x_1, x_2 \right)$
    \begin{tmindent}
      $s \leftarrow 0, x \leftarrow x_2, \beta \leftarrow \left\lfloor n /
      \left( x_{} + 1 \right) \right\rfloor, \varepsilon \leftarrow n
      \tmop{mod} \left( x + 1 \right), \delta \leftarrow \left\lfloor n / x_{}
      \right\rfloor - \beta, \gamma \leftarrow \beta - x_{} \delta$
      
      \tmtextbf{while} $x \geq x_1$ \tmtextbf{do}
      \begin{tmindent}
        $\varepsilon \leftarrow \varepsilon + \gamma$
        
        \tmtextbf{if} $\varepsilon \geq x$ \tmtextbf{then}
        \begin{tmindent}
          $\delta \leftarrow \delta + 1, \gamma \leftarrow \gamma - x,
          \varepsilon \leftarrow \varepsilon - x$
          
          \tmtextbf{if} $\varepsilon \geq x$ \tmtextbf{then}
          \begin{tmindent}
            $\delta \leftarrow \delta + 1, \gamma \leftarrow \gamma - x,
            \varepsilon \leftarrow \varepsilon - x$
          \end{tmindent}
          \begin{tmindent}
            \tmtextbf{if} $\varepsilon \geq x$ \tmtextbf{then} \tmtextbf{exit}
            \tmtextbf{while} \tmtextbf{end} \tmtextbf{if}
          \end{tmindent}
          \tmtextbf{end} \tmtextbf{if}
        \end{tmindent}
        \tmtextbf{else} \tmtextbf{if} $\varepsilon < 0$ \tmtextbf{then}
        \begin{tmindent}
          $\delta \leftarrow \delta - 1, \gamma \leftarrow \gamma + x,
          \varepsilon \leftarrow \varepsilon + x$
        \end{tmindent}
        \tmtextbf{end} \tmtextbf{if}
        
        $\gamma \leftarrow \gamma + 2 \delta, \beta \leftarrow \beta + \delta,
        s \leftarrow s + \beta, x \leftarrow x - 1$
      \end{tmindent}
      \tmtextbf{end} \tmtextbf{while}
      
      $\varepsilon \leftarrow n \tmop{mod} \left( x + 1 \right), \delta
      \leftarrow \left\lfloor n / x_{} \right\rfloor - \beta, \gamma
      \leftarrow \beta - x_{} \delta$
      
      \tmtextbf{while} $x \geq x_1$ \tmtextbf{do}
      \begin{tmindent}
        $\varepsilon \leftarrow \varepsilon + \gamma, \delta_2 \leftarrow
        \left\lfloor \varepsilon / x \right\rfloor, \delta \leftarrow \delta +
        \delta_2, \varepsilon \leftarrow \varepsilon - x \delta_2$
        
        $\gamma \leftarrow \gamma + 2 \delta - x \delta_2, \beta \leftarrow
        \beta + \delta, s \leftarrow s + \beta, x \leftarrow x - 1$
      \end{tmindent}
      \tmtextbf{end} \tmtextbf{while}
      
      \tmtextbf{while} $x \geq x_1$ \tmtextbf{do}
      \begin{tmindent}
        $s \leftarrow s + \left\lfloor n / x \right\rfloor, x \leftarrow x -
        1$
      \end{tmindent}
      \tmtextbf{end} \tmtextbf{while}
      
      \tmtextbf{return} s
    \end{tmindent}
    \tmtextbf{end} \tmtextbf{function}
  \end{tmindent}
\end{tmparsep}{\hspace*{\fill}}{\medskip}

\section{Time and Space Complexity}

Now we present an analysis of the runtime behavior of algorithm.

\begin{theorem}
  The time complexity of algorithm [\ref{algorithm1}] when computing $T \left(
  n \right)$ is $O \left( n^{1 / 3} \right)$ and the space complexity is $O
  \left( \log n \right)$.
\end{theorem}

Before we start, we realize that because $x_{\min} = O \left( n^{1 / 3}
\right)$ and we handle the values of $1 \leq x < x_{\min}$ manually, the
algorithm is at best $O \left( n^{1 / 3} \right)$. In this section we desire
to show that the rest of the computation is at worst $O \left( n^{1 / 3}
\right)$ so that this lower bound holds for the entire computation.

Our first task is to count and size all the top-level regions. We process one
top level region for each integral slope $- a$ from $- 1$ to the slope at
$x_{\min}$. The value for $a$ at each value of $x$ is given by:
\begin{equation}
  a_{} = - \frac{d}{d x} Y \left( x_{} \right) = \frac{n}{x^2}
\end{equation}
and:
\begin{equation}
  X \left( a \right) = \sqrt{\frac{n}{a}}
\end{equation}
Choosing $C_1 = 1$ so that $x_{\min} = \sqrt[3]{2 n}$, then the highest value
of $a$ processed is:
\begin{equation}
  a_{\max} = \frac{n}{x_{\min}^2} = \frac{n^{1 / 3}}{2^{2 / 3}}
\end{equation}
so there are $O \left( n^{1 / 3} \right)$ top level regions.

How big is each top level region? The change in $x$ per unit change in $a$ is
$d x / d a$ and so:
\begin{equation}
  A = - \frac{d}{d a} X \left( a \right) = \frac{n^{1 / 2}}{2 a^{3 / 2}}
\end{equation}
Assume for the moment that the number of total regions visited while
processing a region of size $A$ is:
\[ N \left( A \right) = O \left( A^G \right) \]
noting that the cost of processing a region (excluding the cost of processing
its sub-regions) is $O \left( 1 \right)$ and so the total number of regions is
representative of the total cost.

Now we sum the number of sub-regions processed across all top level region:
\begin{eqnarray}
  N_{\tmop{total}} = \sum^{a_{\max}}_{a = 2} N \left( A \right) & = & O \left(
  \int^{a_{\max}}_1 N \left( A \right) d a \right) \nonumber\\
  & = & O \left( \int_1^{a_{\max}} \left( \frac{n^{1 / 2}}{2 a^{3 / 2}}
  \right)^G d a \right) 
\end{eqnarray}
We can classify three cases depending on the value of $G$ because the outcome
of the integration depends on the final exponent of $a$:
\[ N_{\tmop{total}} = \left\{ \begin{array}{ll}
     O \left( n^{1 / 3} \right) & \tmop{if} G < 2 / 3 ;\\
     O \left( n^{1 / 3} \log n \right) & \tmop{if} G = 2 / 3 ;\\
     O \left( n^{G / 2} \right) & \tmop{if} G > 2 / 3.
   \end{array} \right. \]
(Note that we cannot get below $O \left( n^{1 / 3} \right)$ even if $G = 0$
because we have at least $a_{\max} = O \left( n^{1 / 3} \right)$ top level
regions.)

Now let us analyze the exponent in $N \left( A \right)$. In order to determine
the number of regions encountered in the course of processing a region of size
$A$, we need to analyze the recursion depth. The recursion will terminate when
$w$ or $h$ is unity because by our conditions it is then impossible for the
region to contain any more lattice points. Our next task is to measure the
size of such a region and so we need to know how many $x$ lattice columns that
terminal region represents.

We can use the transformation between $u v$ and $x y$ coordinates given by
(\ref{eq:uv2xy1}) to compute the difference between the $x$ coordinates of
$P_2$ at $\left( 1, 0 \right)$ and $P_1$ at $\left( 0, 1 \right)$, assuming
the smallest case with $w = h = 1$:
\begin{equation}
  \Delta x = x_2 - x_1 + 1 \geq \left( x_0 + 1 \cdot b_2 - 0 \cdot b_1 \right)
  - \left( x_0 + 0 \cdot b_2 - 1 \cdot b_1 \right) + 1 = b_1 + b_2 + 1 > b_1 +
  b_2
\end{equation}
so the size of a terminal region is greater than the sum of the denominators
of the slopes of the two lines that define it.

Each time we recurse into two new regions we add a new extended Farey fraction
that is the mediant of the two slopes for the outer region. As a result, we
perform a partial traversal of a Stern-Brocot tree, doubling the number of
nodes at each level. However, for our current purposes we can ignore the
numerators because we are interested in the sum of denominators. \ Because
regions cannot overlap, this means that the sum of the denominators at the
deepest level of the tree cannot exceed the size of the first region and that
only denominators affect the recursion depth.

Next we need to derive a formula for the sum of the denominators of a partial
Stern-Brocot tree of depth $D$. For example, if the first node $\left( a_1 /
b_1, a_2 / b_2 \right)$ is $\left( 2 / 1, 1 / 1 \right)$, the next two nodes
are $\left( 2 / 1, 3 / 2 \right)$ and $\left( 3 / 2, 1 / 1 \right)$.
Continuing and ignoring numerators we have the following $\left( b_1, b_2
\right)$ tree:

\includegraphics{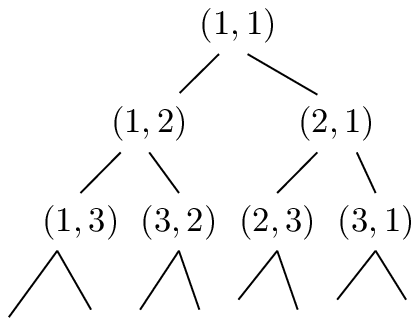}

At each new level we have twice as many nodes and half of the numbers are
duplicated from the previous level and the other half of the numbers are the
sum of numbers of their parent node. Since each parent's sum contributes to
exactly two numbers in the children, the sum of the denominators at each level
is triple the sum of the previous level. So staring with $1 + 1 = 2$ leads to
the sequence $2, 6, 18, 54, \ldots$, and denoting by $\Omega$ the set of
terminal regions, the sum at depth $D$ is therefore
\[ A > \sum_{R : R \in \Omega} b_1 + b_2 = 23^D . \]
Because the number of terminal regions is $\left| \Omega \right| = 2^D$, we
can now place a bound on $\left| \Omega \right|$ in terms of $A$:
\[ \left| \Omega \right| < \left( \frac{A}{2} \right)^{1 / \log_2 3} . \]
Finally, since the total number of regions is $1 + 2 + 4 + \ldots + \left|
\Omega \right| = \sum^D_{i = 1} 2^i$, the number of regions as a function of
the size $A$ is
\begin{equation}
  N \left( A \right) = 2 \left| \Omega \right|_{} - 1 = O \left( A^{1 / \log_2
  3} \right)
\end{equation}
and therefore $G = 1 / \log_2 3$.

Since $1 / \log_2 3 \approx 0.63$, this means that $G < 2 / 3$ and the proof
that the overall time complexity of the algorithm is $O \left( n^{1 / 3}
\right)$ is complete.

The space complexity is simply our recursion depth which can be at most $O
\left( \log n \right)$.

\section{Higher-Order Divisor Sums}

The two-dimensional hyperbola and the functions $\tau \left( n \right)$ and $T
\left( n \right)$ can be generalized to higher dimensions. Using this notation
$\tau \left( n \right) = \tau_2 \left( n \right)$ and $T \left( n \right) =
T_2 \left( n \right)$. Then the divisor sum $T_3 \left( n \right)$, the
summatory function for $\tau_3 \left( x \right) = \sum_{abc = x} 1$, can be
computed by summing under the three-dimensional hyperbola
\[ T_3 \left( n \right) = \sum_{x, y, z : xyz \leq n} 1 = \sum_{z = 1}^n
   \sum_{x = 1}^n \left\lfloor \frac{n}{xz} \right\rfloor = \sum_{z = 1}^n T
   \left( \left\lfloor \frac{n}{z} \right\rfloor \right) . \]
Again using the symmetry of this hyperbola we can restrict the outer summation
to $\sqrt[3]{n}$ by counting nested ``shells'', and avoiding double and triple
counting, we get
\begin{eqnarray*}
  T_3 \left( n \right) & = & \sum_{z = 1}^{\left\lfloor \sqrt[3]{n}
  \right\rfloor} \left[ 3 \left( 2 \sum^{\left\lfloor \sqrt{\frac{n}{z}}
  \right\rfloor}_{x = z + 1} \left( \left\lfloor \frac{n / z}{x} \right\rfloor
  - z \right) - \left( \left\lfloor \sqrt{\frac{n}{z}} \right\rfloor - z
  \right)^2 + \left( \left\lfloor \frac{n}{z^2} \right\rfloor - z \right)
  \right) + 1 \right]\\
  & = & \sum_{z = 1}^{\left\lfloor \sqrt[3]{n} \right\rfloor} \left[ 3 \left(
  2 \sum^{\left\lfloor \sqrt{\frac{n}{z}} \right\rfloor}_{x = z + 1}
  \left\lfloor \frac{n / z}{x} \right\rfloor - 2 z \left( \sqrt{\frac{n}{z}} -
  z \right) - \left( \left\lfloor \sqrt{\frac{n}{z}} \right\rfloor^2 - 2 z
  \sqrt{\frac{n}{z}} + z^2 \right) + \left\lfloor \frac{n}{z^2} \right\rfloor
  - z \right) + 1 \right]\\
  & = & \sum_{z = 1}^{\left\lfloor \sqrt[3]{n} \right\rfloor} \left[ 3 \left(
  2 S \left( \left\lfloor \frac{n}{z} \right\rfloor, z + 1, \left\lfloor
  \sqrt{\frac{n}{z}} \right\rfloor \right) - \left\lfloor \sqrt{\frac{n}{z}}
  \right\rfloor^2 + \left\lfloor \frac{n}{z^2} \right\rfloor + z^2 - z \right)
  + 1 \right]\\
  & = & 3 \sum_{z = 1}^{\left\lfloor \sqrt[3]{n} \right\rfloor} \left( 2 S
  \left( \left\lfloor \frac{n}{z} \right\rfloor, z + 1, \left\lfloor
  \sqrt{\frac{n}{z}} \right\rfloor \right) - \left\lfloor \sqrt{\frac{n}{z}}
  \right\rfloor^2 + \left\lfloor \frac{n}{z^2} \right\rfloor \right) +
  \left\lfloor \sqrt[3]{n} \right\rfloor^3
\end{eqnarray*}
where in the last step we use the identity $\sum_{z = 1}^k 3 \left( z^2 - z
\right) + 1 = 3 \left( k \left( k + 1 \right)  \left( 2 k + 1 \right) / 6 - k
\left( k + 1 \right) / 2 \right) + k = k^3$. Since $S \left( n, x_1,
\left\lfloor \sqrt{n} \right\rfloor \right)$ is a partial result in the
calculation of $T \left( n \right)$, it is also has $O \left( n^{1 / 3}
\right)$ time complexity when using Algorithm [\ref{algorithm1}]. \ As a
result, we can calculate $T_3 \left( n \right)$ in
\[ \sum^{\left\lfloor \sqrt[3]{n} \right\rfloor}_{z = 1} O \left( \left\lfloor
   \frac{n}{z} \right\rfloor^{1 / 3} \right) = O \left( \int_1^{n^{1 / 3}}
   \frac{n^{1 / 3}}{z^{1 / 3}}^{} d z \right) = O \left( n^{5 / 9} \right), \]
a modest improvement over $O \left( n^{2 / 3} \right)$ using a direct double
summation. \ Similar derivations give $O \left( n^{2 / 3} \right)$ for $T_4
\left( n \right)$ and $O \left( n^{11 / 15} \right)$ for $T_5 \left( n
\right)$ or $O \left( n^{1 - 4 / \left( 3 k \right)} \right)$ for $T_k \left(
n \right)$ in general.

\section{Remarks}.

It would be possible to simplify the algorithm somewhat by removing the
distinction between top level regions and region processing itself by starting
with the region defined by $\left( 1 / 0, 1 / 1 \right)$. The reason for the
current assymetry is two-fold. First, some of the solutions to the equations
are degenerate when $a_i b_i = 0$ and would require special handling anyway.
Second, and perhaps more importantly, we can also capitalize on the simpler $x
y$ coordinate system where possible.

The two major sections of the algorithm, $S_1$ and $S_4$, are easily
parallelizable. The section $S_1$ can divide summation batches to different
processors. The section $S_4$ can be revised to use a work queue of regions
instead of recursion. During region processing, one region can be enqueued and
the other processed iteratively. Available processors can dequeue regions that
need to be processed.

In fact it turns out that the $S \left( n / z, z + 1 \left\lfloor \sqrt{n / z}
\right\rfloor \right)$ terms in the $T_3 \left( n \right)$ summation skip over
the problematic first $O \left( \left( n / z \right)^{1 / 3} \right)$ columns
by the time $z$ reaches $n^{1 / 4}$ and then start eroding away the smallest
regions as $z$ approaches $n^{1 / 3}$. \ Modifying the method slightly and
then computing the time complexity of these two portions separately and
allowing $a_{\max}$ to decline appropriately we would achieve $O \left( n^{1 /
2} \log n \right)$ for $T_3 \left( n \right)$ if we could prove that $G = 1 /
2$. \ In any case, using $G = 1 / \log_2 3$ at least gives us $O \left( n^{5 /
9 - c + \epsilon} \right)$ for some $c > 0$.

\section{Related Work}

In [\ref{bib:Gal00}], Galway presents an improved sieving algorithm that also
features region decomposition based on extended Farey fractions as well as
coordinate transformation. In [\ref{bib:Tao11}], applications for the divisor
summatory are function presented including computing the parity of $\pi \left(
x \right)$, the prime counting function, as well as a sketch for a different
$O \left( n^{1 / 3} \right)$ algorithm. In [\ref{bib:Sil12}], the parity of
the prime counting function is studied more closely and several related
algorithms are developed.


\begin{thebibliography}{}
  
  
  \bibitem{Vor03}\label{bib:Vor03}Georges Vorono\"{\i}, \tmtextit{Sur un
  probl\`eme du calcul des fonctions asymptotiques}, J. Reine Angew. Math.
  \tmtextbf{126} (1903), 241-282.
  
  \bibitem{Bres77}\label{bib:Bres77}Jack Bresenham, \tmtextit{A linear
  algorithm for incremental digital display of circular arcs}, Communications
  of the ACM \tmtextbf{20} (1977), 100-106.
  
  \bibitem{Gal00}\label{bib:Gal00}William F. Galway, \tmtextit{Dissecting a
  Sieve to Cut Its Need for Space}, In Proceedings of ANTS. (2000), 297-312.
  
  \bibitem{Tao11}\label{bib:Tao11}Terence Tao, Ernest Croot III, and Harald
  Helfgott. \tmtextit{Deterministic methods to find primes}. Mathematics of
  Computation, 2011. Published electronically on August 23, 2011.
  
  \bibitem{Sil12}\label{bib:Sil12}Tom\'as Oliveira e Silva,
  \tmtextit{Efficient Computation of the Parity of the Prime Counting
  Function}, in preparation.
\end{thebibliography}
\end{document}